\documentclass[12pt,twoside]{amsart}
\usepackage[dvips]{graphics}
\usepackage[all]{xy}
\usepackage{a4,amsmath,amssymb,amsfonts,amscd,mathrsfs}
\reversemarginpar
\newcounter{Abschnitt}[section]

\newtheorem{theorem}[subsection]{Theorem}
\newtheorem{lemma}[subsection]{Lemma}
\newtheorem{corollary}[subsection]{Corollary}
\newtheorem{proposition}[subsection]{Proposition}

\renewcommand{\Im}{\operatorname{Im}}

\newcommand{\Orth}{{\rm O}}

\renewcommand{\Re}{\operatorname{Re}}

\newcommand{\dee}{\partial}

\renewcommand{\phi}{\varphi}

\newcommand{\cdop}{{\mathbb C}}
\newcommand{\edop}{{\mathbb E}}

\newcommand{\ndop}{{\mathbb N}}

\newcommand{\qdop}{{\mathbb Q}}
\newcommand{\rdop}{{\mathbb R}}
\newcommand{\zdop}{{\mathbb Z}}

\author{Juan Marcos Cervi\~no}
\address{FB Mathematik, University Duisburg-Essen,
45117 Essen, Germany}
\email{juan.cervino@uni-due.de}
\author{Georg Hein}
\address{FB Mathematik, University Duisburg-Essen,
45117 Essen, Germany}
\email{georg.hein@uni-due.de}
\date{June 5, 2009}
\begin{document}
\title{lattice invariants from the heat kernel}
\maketitle
\begin{abstract}
We derive lattice invariants from the heat flux of a lattice.
Using systems of harmonic polynomials, we obtain sums of products
of spherical theta functions which give new invariants of integer
lattices which are modular forms. In particular, we show that the
modular forms $\Theta_{nn,\Lambda}$
depend only from lengths and angles in the lattice.
\end{abstract}
\section{Introduction}
For a lattice $\Lambda \subset \edop^n$ one of the most important
invariants is its theta function $\Theta_\Lambda$.
This function does not determine the lattice.
We want to describe and illustrate a way to derive more lattice
invariants from the heat flux of the lattice. The first idea is to apply
the heat flux to the lattice $\Lambda$ considered as distribution on
functions on $\edop^n$. In doing so, we obtain a function $f_t
=f_{t,\Lambda}: \rdop^+ \times \edop^n \to \rdop$. Secondly, we try to
deduce invariants from $f_t$ which do not depend on the embedding of the
lattice. The restriction
of $f$ to $\rdop^+ \times \{ 0 \}$ is a first such invariant. However it
gives essentially the theta function of the lattice (see
\ref{first-ex}). An obvious way to obtain invariants is to integrate
products homogeneous parts of $f$ over the unit spehere $S^{n-1} \subset
\edop^n$. These invariants are denoted by $c_{k_1,\dots,k_m}$ (see
\ref{defck1km}).
In order to obtain modular forms, we study polynomial differentials of
$f$, in particular for harmonic polynomials (see Proposition
\ref{harm-diff}).  Whenever we find an invariant sum of such products
(which we call a harmonic datum), then we obtain an invariant modular
form (see Proposition \ref{data-theta}). Thus, we need a construction
for harmonic data. We illustrate how the invariants $c_{k_1,\dots,k_m}$
can be used to construct harmonic data in Section \ref{thetann-2} where
we construct the harmonic data $p_{11}$ and $p_{22}$ for lattices in
$\edop^2$. Seeing the results we obtain from these data (cf. Corollaries
\ref{theta1} and \ref{theta2}) it is not hard to show that for a lattice
$\Lambda \subset \edop^2$ of level $N$ the holomorphic function
\[ \Theta_{nn,\Lambda}(\tau) = \sum_{(\gamma,\delta) \in \Lambda \times
\Lambda} \left(
\cos(2n\measuredangle
(\gamma,\delta ))||\gamma||^{2n}||\delta ||^{2n}
\right)q^{||\gamma||^2+||\delta||^2}
\quad \mbox{ with } q=\exp(2 \pi i \tau) \]
is a modular form with integer coefficients of weight $4n+2$ and of
level $N$ (see Theorem \ref{theta3} for more details). In the final
Section \ref{theta11-n} we start to generalize the modular forms
$\Theta_{mm,\Lambda}$  to lattices in $n$ dimensional euclidean space
$\edop^n$. In particular we deduce in Theorem \ref{theta11} a concrete
formula for $\Theta_{11,\Lambda}$ which we apply to two isospectral
lattices $\Lambda_1$ and $\Lambda_2$
in $\edop^4$ found by A.~Schiemann in \cite{Sch}. An elementary
calculation (see Proposition \ref{exam-schie}) shows that
$\Theta_{11,\Lambda_1} \ne \Theta_{11,\Lambda_2}$. Thus, the modular
lattice invariant $\Theta_{11,\Lambda}$ can distinguish at least two
isospectral lattices.\\
{\em Notation: } By a lattice we mean a free $\zdop$ module $\Lambda
\cong \zdop^n$ with a positive definite quadratic form
$q:\Lambda \to \rdop$. A lattice possesses isometric
embeddings $\Lambda \to \edop^n$. Two such embeddings differ by an
element of the orthogonal group $\Orth(n)$. The number
$q(\gamma)$ is the length of $\gamma \in \Lambda$. The norm of
$\gamma \in \Lambda$ is defined to be
$||\gamma||:=\sqrt{q(\gamma)}$.

\newpage
\section{The function $f_t$}\label{HEAT}
\subsection{The heat flux of a lattice}\label{heat}
We consider a cocompact lattice $\Lambda \subset \edop^n$ in the $n$
dimensional euclidean space. We consider two lattices $\Lambda,\Lambda'
\subset \edop^n$ to be isomorphic when there exists an element $\phi \in
\Orth(n)$ which transforms $\Lambda$ into $\Lambda'$, in short
$\Lambda'=\phi(\Lambda)$.

We regard $\Lambda \in D^n(\edop^n)$ as a distribution on the
Schwartz space
$S(\edop^n)$ by
\[ \Lambda(f)= \sum_{\gamma \in \Lambda} f(\gamma). \]
Smooth $n$ forms $\omega$ on $\edop^n$ can be considered as currents
$[\omega] \in D^n(\edop^n)$ by $[\omega]f:= \int_{\edop^n} \omega \cdot
f$. A current $[\omega]$ associated to a smooth form $\omega$ is called
a smooth current.
We want to subdue the singular distribution $\Lambda$ to the heat flux
to obtain a smooth current $(\Lambda)_t$ for any positive $t \in
\rdop$. We follow the approach in \cite[Section 2]{Hei}:
\[ \mbox{The Laplace operator is given by:} \qquad
 \Delta = -\sum_{i=1}^n \frac{\dee^2}{\dee x_i^2} \, . \]

For any $k$-form $f \in A^k(\edop^n)$ we have the heat flux of $f$:
\[ P_t(f) \in A^k(X) \mbox{ for all } t \in \rdop_+  \quad
 \left(\frac{\dee}{\dee t}+\Delta \right) P_t(f)=0  \quad
 \lim_{t \to 0} P_t(f) = f.
\]
The heat flux of a distribution  $T \in D^n(\edop^n)$ is defined by
\[ P_t(T)(f):= T(P_t(f)).\]
If the distribution is smooth, then we have the equality:
\[P_t([\omega]) = [ P_t(\omega)] . \]
That is the extension to distributions is compatible with the embedding
$A^n(\edop^n) \hookrightarrow D^n(\edop^n)$ and the heat flux on smooth
forms. Furthermore, for $T=\Lambda$ the heat flux smoothes the
distribution, that is
\[ P_t(\Lambda) = [\omega_{t}] \quad
\mbox{for some smooth } \omega_{t} \in A^n(\edop^n) .\]
Using the explicit form of the heat kernel for a point (cf.
\cite[Chapter 2]{BGV}) we can write down the heat flux of $\Lambda$
directly as
$\omega_{t} = f_{t} d \mu$, with 
$d\mu =dx_1 \wedge \dots dx_n$ the volume form on $\edop^n$, and
\[ f_{t}(x) = {(4\pi t)^\frac{-n}{2}}\sum_{\gamma \in \Lambda}
\exp \left( \frac{-||x-\gamma||^2}{4t} \right) \, . \]

\subsection{$f_t$ determines the lattice $\Lambda \subset \edop^n$}
Indeed we have on the level of distributions $\lim_{t \to +0} [f_t d
\mu] = \Lambda $.
However, the function $f_t$ 
depends on the embedding of $\Lambda$.  We try to deduce more informations by considering the values
of derivatives of $f_t$ along $\rdop^+ \times \{ 0 \}$.

\subsection{Polynomial derivations}\label{pairing}
We need the subring $A = \rdop[x_1,\dots,x_n] \subset A^0(\edop^n)$ of
polynomials in the smooth functions. We define a pairing on
\[ A \times A^0(\edop^n) \to \rdop \quad \langle P , f \rangle =
P\left(\frac{\dee}{\dee x_1},\dots \frac{\dee}{\dee x_n}\right)f |_0 \,.\]
Note that the pairing is also well defined when $P$ is smooth and $f$
is a polynomial. We have the following
properties of that pairing:\\
\begin{tabular}{lp{14cm}}
(i) & The pairing is bilinear.\\
(ii) & The pairing is symmetric, in the sense that $\langle P , f
\rangle = \langle f , P \rangle$ whenever one side is defined.\\
(iii) & the restriction of $\langle \, , \, \rangle$ to $A
\times A$ is positive definite.\\
(iv) & The monomials form a orthogonal (not orthonormal!) basis.\\
(v) & For two polynomials $P,Q \in A$ we have $\langle P \cdot Q ,f
\rangle= \langle P , Q\left(\frac{\dee}{\dee x_1},\dots \frac{\dee}{\dee
x_n}\right)f \rangle$.
\end{tabular}
For any $P \in A$ we obtain with $\langle P , f_t \rangle$ a smooth
function $\langle P , f_t \rangle: \rdop^+ \to \rdop$.

\begin{lemma}\label{delta=dt}
For the $\Orth(n)$-invariant polynomial $r^2:=\sum_{i=1}^n x_1^2$ we
have
\[ \langle r^{2k} , f_t \rangle = \frac{\dee^k}{\dee t^k}f_t(0) \, .\]
\end{lemma}
\begin{proof} Having in mind, that $f_t$ is the heat flux of $\Lambda$,
that is $(\Delta+\frac{\dee}{\dee t} ) f_t =0$, we get
\[ \langle r^2 , f_t \rangle = \langle 1 , r^2\left(\frac{\dee}{\dee
x_1},\dots \frac{\dee}{\dee x_n}\right) f_t \rangle =
\langle 1 , -\Delta(f_t) \rangle
= \langle 1 , \frac{\dee}{\dee t}f_t \rangle = \frac{\dee}{\dee t}f_t(0)
. \]
Now the assertion of the lemma follows from \ref{pairing} (v).
\end{proof}
So evaluation of $\Orth(n)$ invariant polynomials $P=\sum_i a_i r^{2i}$
gives only expressions of type $\sum_i a_i \frac{\dee^i}{\dee t^i}
f_t(0)$.
After having seen what $\langle r^{2k},f_t \rangle$ is, we study the
pairing $\langle h,f_t \rangle$ for homogeneous harmonic polynomials.
We have the following result:

\begin{proposition}\label{harm-diff}
If $h \in A$ is an homogeneous harmonic polynomial of degree $d$, then
\[ (2t)^d\langle h,f_t \rangle 
=(4 \pi t)^{\frac{-n}{2}}  \sum_{\gamma
\in \Lambda} h(\gamma) \exp \left( \frac{-||\gamma||^2}{4t}\right) .\]
\end{proposition}

In order to prove this proposition we need an auxiliary result:

\begin{lemma}\label{ha-di-le}
Let $h$ be homogeneous polynomial of degree $d$, $a \in \cdop^*$ a number
and $f=f(x_1,\dots,x_n) = \exp\left( \frac{a}{2} \cdot \sum_{i=1}^n x_i^2
\right)$. Then we have an equality
\[ h\left( \frac{\dee}{\dee x_1}, \dots, \frac{\dee}{\dee x_1}\right) f = 
a^d \left( \sum_{k \geq 0} \left(\frac{-1}{2a} \right)^k \frac{1}{k!} \, \Delta^k(h)
\right) \cdot f .\]
\end{lemma}
{\em Proof.}
We have to prove the lemma only for $h$ a monomial. We proceed by
induction on the degree of $h$. For $\deg(h)=0$, the statement is
obvious. Suppose the assertion holds for $h'$ of degree $d-1$.
We set $h=x_ih'$. An direct calculation shows that
\[ \Delta(h) = x_i \Delta(h') -2 \frac{\dee}{\dee x_i}
h' .  \]
Since $\Delta$ and $\frac{\dee}{\dee x_i}$ commute we deduce from that
formula inductively that
\begin{equation}
\Delta^k(h) = x_i \Delta^k(h') -2k \frac{\dee}{\dee x_i}\Delta^{k-1}(h')
\end{equation}
holds for all integers $k \geq 0$. 
Now we compute using the induction hypothesis, the Leibniz rule, and
equation (1):
\[ \begin{array}{rclr}
h\left( \frac{\dee}{\dee x_1}, \dots, \frac{\dee}{\dee x_1}\right) f
&=& \frac{\dee}{\dee x_i} \left(  a^{d-1} \left( \sum_{k \geq 0}
\left(\frac{-1}{2a} \right)^k \frac{1}{k!} \, \Delta^k(h')
\right) \cdot f \right)\\
&= & a^d \left( \sum_{k \geq 0} \left(\frac{-1}{2a} \right)^k
\frac{1}{k!} (x_i\Delta^k(h') - 2k \frac{\dee}{\dee x_i}
\Delta^{k-1}(h')) \right)\cdot f\\
&=& a^d \left( \sum_{k \geq 0} \left(\frac{-1}{2a} \right)^k \frac{1}{k!}
\, \Delta^k(h)
\right) \cdot f \, . & \Box \\
\end{array} \]

{\em Proof of Proposition \ref{harm-diff}.}
Setting $a=\frac{-1}{2t}$ in Lemma \ref{ha-di-le} we obtain that
\[ h\left( \frac{\dee}{\dee x_1}, \dots, \frac{\dee}{\dee x_1}\right)
\exp \left( \frac{-\sum_{i=1}^n x_i^2 }{4t} \right)  = \left( \frac{-1}{2t} \right)^d
h(x_1,\dots,x_n) \exp \left(\frac{-\sum_{i=1}^n x_i^2 }{4t}  \right) \, .\]
Setting $x_i = \tilde x_i- \gamma_i$ for some $\gamma =(\gamma_1,\dots
,\gamma_n ) \in \Lambda$ this reads
\[ h\left( \frac{\dee}{\dee x_1}, \dots, \frac{\dee}{\dee x_1}\right)
\exp \left( \frac{-||\gamma-x||^2}{4t} \right) = \left( \frac{-1}{2t}
\right)^d h(x_1-\gamma_1,\dots,x_n-\gamma_n) \exp \left(
\frac{-||\gamma-x||^2}{4t} \right)   \]
Multiplying this equality with $(4\pi t)^{\frac{-n}{2}}$ and summing up
over all $\gamma \in \Lambda$ yields the statement of the proposition
when specializing to $x=0$.
\hfill $\Box$

\subsection{Definition of the $c_{k_1,k_2,\ldots,k_m} $}\label{defck1km}
In order to find an invariant we first decompose the function $f_t$ by
in its homogeneous parts $f_t = f_0 +f_1 + \dots$ with
\[ f_k = \sum _{I \subset \ndop^n, |I|=2k} a_I \frac{x^I}{I!} \mbox{
where for } I=(i_1,i_2,\dots,i_n) \,\,\,
I!:=\prod_{m=1}^n i_m!  \mbox{ ,  } x^I:=\prod_{m=1}^n x_m^{i_m}
\,\mbox{ , and }\]
\[ a_I:= \langle x^I , f_t \rangle = \frac{\dee^{i_1}}{\dee x_1^{i_1}}
\frac{\dee^{i_2}}{\dee x_2^{i_2}} \cdots \frac{\dee^{i_n}}{\dee
x_n^{i_n}}f_t|_{x=0} \, \] 
These polynomials are of course not invariant but the integrals of their products
\[ c_{k_1,k_2,\ldots,k_m} : = 
 \int_{S^{n-1}} f_{k_1} \cdot f_{k_2} \dots f_{k_m} d \bar \mu \]
are. Here $d \bar \mu$ denotes the normalized $\Orth(n)$ invariant
measure on $S^{n-1}$ such that $\int_{S^{n-1}} d \bar  \mu=1$.
We restrict to the parts of even degree because the parts of odd degree
vanish for $f_t$ being an even function.
We will see that the $c_{k_1,k_2,\ldots,k_m}$ can be expressed as
products of expressions of type $\langle h, f_t \rangle$.
The decomposition of a polynomial $f \in A$ into a sum
$f=h_n+r^2h_{n-2} + r^4 h_{n-4} + \dots$ with $h_i$ harmonic will
allow applications of the following:

\subsection{Spherical theta function principle}
We want to compute these invariants and relate them to modular
forms for integral lattices. Here we call a lattice $\Lambda \subset
\edop^n$ integral when $||\gamma||^2 \in \ndop$ for all $\gamma \in
\Lambda$. Such a lattice has a level $N \in \ndop$ and a discriminant
$D$ (cf. \cite[Chapter 3.1]{Zag}). For an homogeneous harmonic function $h$
on $\edop^n$ of degree $d$ we have the spherical theta function
$\Theta_{h,\Lambda}$ defined by
\[ \Theta_{h,\Lambda}(\tau) = \sum_{ \gamma \in \Lambda}
h(\gamma)q^{||\gamma||^2} \, ,\qquad \qquad \qquad q=\exp(2 \pi i \tau). \]
These spherical theta functions are modular forms of
weight $\frac{n}{2}+d$, level $N$, and
character $\left(\frac{D}{} \right)$ (see again \cite[Chapter
3.2]{Zag}).\\

{\em Definition.}
An harmonic datum $p=((h_{ij})_{i=1,..,m \,\, j=1,\dots,k})$ for lattices
in $\edop^n$ consists of harmonic homogeneous polynomials $h_{ij}$ of
degree $d_i$ on $\edop^n$,
which define a map
\[ p :\{ \mbox{lattices } \Lambda \subset
\edop^n\} \to  \{ \mbox{smooth functions }  \rdop^+ \to \cdop \} \qquad \Lambda
\mapsto p(\Lambda) \]
\[\mbox{by }  p(\Lambda):= \sum_{j=1}^k \prod_{i=1}^m \langle h_{ij},f_{t,\Lambda}
\rangle 
\qquad \mbox{where} \qquad f_{t,\Lambda}(x) = {(4\pi t)^\frac{-n}{2}}
\sum_{\gamma \in \Lambda}
\exp \left( \frac{-||x-\gamma||^2}{4t} \right)  ,
\]
which satisfies $p(\Lambda) = p(\phi(\Lambda))$ for all $\phi \in
\Orth(n)$.
The sum $d=\sum_{i=1}^m d_i$ is called the degree of the harmonic
datum. The harmonic datum is called even (resp.~odd) when $m$ is even
(resp.~odd).\\

The connection between harmonic data and invariant modular forms
is given the by the following proposition.

\begin{proposition}\label{data-theta} 
Suppose we have a harmonic datum $(p,(h_{ij})_{i=1,..,m \,\,
j=1,\dots,k})$ for lattices in $\edop^n$ of degree $d$. Then for any
integral lattice $\Lambda \subset \edop^n$ of level $N$ the modular form 
\[\Theta_{p,\Lambda}(\tau)
:= \sum_{j=1}^k \prod_{i=1}^m \Theta_{h_{ik},\Lambda} \]
is invariant under $\Orth(n)$.
$\Theta_{p,\Lambda}(\tau)$ is a modular form of weight $\frac{nm}{2}+d$,
and of level $N$.
If the harmonic system is odd,
then $\Theta_{p,\Lambda}(\tau)$ has character
$\left( \frac{D}{} \right)$ with $D$ the discriminant of the lattice.
If it is even, then $\Theta_{p,\Lambda}(\tau)$ is a
modular form for the trivial character.
\end{proposition}
\begin{proof}
First we remark that $\Theta_{p,\Lambda}$ being a sum of products of
modular forms of level $N$ is a modular form of level $N$ and the stated
weight. Likewise we see that it is a modular form of character
$\left(\frac{D}{} \right)^m$.
Thus, it remains to show that $\Theta_{p,\Lambda}$ is
invariant under $\Orth(n)$.\\
We start with calculating $p(\Lambda)(t)$ using Proposition
\ref{harm-diff}:
\[ p(\Lambda)(t) = \sum_{j=1}^k \prod_{i=1}^m \langle
h_{ij},f_{t,\Lambda} \rangle =\sum_{j=1}^k \prod_{i=1}^m (2t)^{-d_i} (4
\pi t)^{\frac{-n}{2}} \sum_{\gamma \in \Lambda} h_{ij}(\gamma) \exp\left(
\frac{-||y||^2}{4t} \right)
 \]
Since the $d_i$ sum up to $d$ we get
\[ (2t)^d (4 \pi t)^{\frac{mn}{2}} p(\Lambda)(t) = \sum_{j=1}^k
\prod_{i=1}^m  \sum_{\gamma \in \Lambda} h_{ij}(\gamma) \exp\left(
\frac{-||y||^2}{4t} \right)
\]
is a lattice invariant. However, this is the value of
$\Theta_{p,\Lambda} \left(\frac{i}{8 \pi t} \right)$, so by the identity
theorem for holomorphic functions the values of the modular form
$\Theta_{p,\Lambda} $ along the imaginary axis determine this function.
\end{proof}

\subsection{Example: The theta series}\label{first-ex}
The function $f_t(0)= {(4\pi t)^\frac{-n}{2}} \sum_{\gamma \in \Lambda}
\exp(-||\gamma||^2/4t)$ is a lattice invariant.
The theta function of the lattice $\Lambda$ is
\[ \Theta_\Lambda(\tau) := \sum_{\gamma \in \Lambda}\exp(2 \pi i \tau
||\gamma||^2) .\]
So up to the scaling factor $(4 \pi t) ^{\frac{-n}{2}}$ the theta
function $\Theta_\Lambda \left(\frac{i}{8 \pi t} \right)$ gives
$f_t(0)$:
\[ f_t(0) = (4 \pi t) ^{\frac{-n}{2}} \Theta_\Lambda \left(\frac{i}{8 \pi
t} \right) .\]
Vice versa, from $f_t(0)$ we can extract the values of $\Theta_\Lambda$
along the imaginary line, which determines $\Theta_\Lambda$ by the identity
theorem for holomorphic functions.

Setting $h_{11} \equiv 1$ we obtain a
harmonic datum $(f_t(0),1)$ of degree zero with $m=k=1$.
This harmonic datum yields the theta function of the lattice up to the
factor $(4\pi t)^\frac{n}{2}$.
We show in the sequel that we can build more harmonic data from our
lattice invariants $c_{k_1,\dots,k_r}$.

\section{Lattices in $\edop^2$}\label{thetann-2}
\subsection{Preparations}
Let $\Lambda \subset \edop^2$ be a fixed lattice.
From Corollary \ref{int-pol} we deduce the following table of integrals
which will be useful for the coming calculations:

\[\begin{array}{rclcrclcrcl}
\int_{S^{1}}x_{0}^{2}  d \bar \mu & = & \frac{1}{2}& \quad &
\int_{S^{1}}x_{0}^{2} x_{1}^{2}  d \bar \mu & = & \frac{1}{8} & \quad &
\int_{S^{1}}x_{0}^{4}  d \bar \mu & = & \frac{3}{8} \\ \\
\int_{S^{1}}x_{0}^{4} x_{1}^{2}  d \bar \mu & = & \frac{1}{16} & \quad &
\int_{S^{1}}x_{0}^{6}  d \bar \mu & = & \frac{5}{16} & \quad &
\int_{S^{1}}x_{0}^{4} x_{1}^{4}  d \bar \mu & = & \frac{3}{128}\\ \\
\int_{S^{1}}x_{0}^{6} x_{1}^{2}  d \bar \mu & = & \frac{5}{128} & \quad &
\int_{S^{1}}x_{0}^{8}  d \bar \mu & = & \frac{35}{128}\\
\end{array}\]
By definition $c_{0} =f_t(0) =\frac{\Theta_\Lambda(\frac{i}{8\pi t})}{4\pi t}$, so we
see that $c_{ 0}$ is essentially the theta function of the lattice.
Next we compute
\[ 4 c_{1} = 4\int_{S_1} \left(a_{20}\frac{x^2}{2}
+ a_{11}xy + a_{02}\frac{y^2}{2} \right) d\bar \mu = a_{20}+a_{02} \,. \]  
Thus, we obtain $c_{1} = \frac{1}{4} \frac{\dee}{\dee t}c_{0}$.
This is a general pattern:
$c_2=\frac{1}{64} \frac{\dee^2}{\dee t^2}c_{0}$,
$c_3=\frac{1}{2304} \frac{\dee^3}{\dee t^3}c_{0}$,
and in general
$c_n=\frac{1}{4^n (n!)^2} \frac{\dee^n}{\dee t^n}c_{0}$.

\subsection{The harmonic datum $p_{11}$ for lattices in
$\edop^2$}\label{p112}
We start with the computation of the lattice invariant $c_{11}$:
\[ \begin{array}{rcl}
32 c_{1,1}&  = & \displaystyle 32
\int_{S^1} \left(a_{20}\frac{x^2}{2} + a_{11}xy + a_{02}\frac{y^2}{2}
\right)^2 d \bar\mu \\ \\
&=& 3a_{20}^2+3a_{02}^2+ 4 a_{11}^2+2a_{20}a_{02}\\
&=& 2(a_{20}+a_{02})^2 + 4 a_{11}^2+(a_{20}-a_{02})^2\\
\end{array}\]\\
The function $p_{11}=32c_{11}-32c_1^2$ is obviously a lattice invariant,
and gives a harmonic datum:
\[ p_{11} = 4\langle x_0x_1, f_t \rangle ^2
+ \langle x_0^2-x_1^2, f_t \rangle ^2 \]
This harmonic datum yields by Proposition \ref{data-theta}:

\begin{corollary}\label{theta1}
For any integral lattice $\Lambda \subset \edop^2$ of level $N$
the modular form
\[ \Theta_{11,\Lambda} =
4\Theta_{x_0x_1,\Lambda}^2+\Theta_{x_0^2-x_1^2,\Lambda}^2 \]
of weight $6$ and level $N$
is invariant under the $\Orth(2)$ action on the embedding $\Lambda \to
\edop^2$. If we write
$\Theta_{11,\Lambda}(\tau) =\sum_{n \in \ndop} a_n q^n$ with
$q =\exp(2 \pi i \tau)$, then the $a_n$ are given by
\[a_n=
\sum_{\tiny \begin{array}{c}
(\gamma,\delta) \in \Lambda \times \Lambda\\
||\gamma||^2+||\delta ||^2=n\\
\end{array}} \cos(2\measuredangle (\gamma,\delta ))||\gamma||^2||\delta ||^2. \]
\end{corollary}
\begin{proof}
The only thing which remains to be shown is the formula for the $a_n$.
This is the result of a straightforward calculation, where we expand
$4\Theta_{x_0x_1,\Lambda}^2$ and $\Theta_{x_0^2-x_1^2,\Lambda}^2$ as sums
over $\Lambda \times \Lambda$, and check that the summand of
$q^{||\gamma||^2+||\delta ||^2}$ is $\cos(2\measuredangle
(\gamma,\delta ))||\gamma||^2||\delta ||^2$ for each pair
$(\gamma,\delta ) \in \Lambda \times \Lambda$.
See the proof of Theorem \ref{theta3} for this calculation.
\end{proof}

\subsection{The harmonic datum $p_{22}$ for lattices in $\edop^2$}
We start with a decomposition of $c_{22}$ into three summands
\[ \begin{array}{rcl}
73728 c_{2,2}&=&
6a_{{40}}a_{04}+96a_{{31}}a_{{13}}+60a_{{22}}a_{{04}}+60a_{{40}}a_{{22}}+108{a_{{22}}}^{2}+\\
&&  \qquad \qquad \qquad \qquad \qquad+ 80{a_{{13}}}^{2}+35{a_{40}}
^{2}+35{a_{{04}}}^{2}+80{a_{{31}}}^{2} \\
&=& p_{22} + s_1 +s_2\\
\mbox{with} \qquad p_{22}&
=&(a_{40}-6a_{22}+a_{04})^2 + 16(a_{31}-a_{13})^2\\
s_1&=& 18(a_{40}+2a_{22}+a_{04})^2\\
s_2&=&16((a_{20}+a_{02})\circ(a_{20}-a_{02}))^2+
64((a_{20}+a_{02})\circ a_{11})^2 \\
\end{array}\] \\
If we can show that $s_1$ and $s_2$ are invariant under the $\Orth(2)$
action,
then it follows that $p_{22}$ is also invariant under that action.
However, we see that $s_1=73728 c_2^2$, which is an invariant.
Next we see that
\[s_2=16\left( 4\langle x_0x_1 , \frac{\dee}{\dee t} f_t
\rangle^2 + \langle x_0^2-x_1^2 , \frac{\dee}{\dee t} f_t
\rangle^2\right) \, . \]
As we have see before Corollary \ref{theta1} this is also 
invariant under the $\Orth(2)$ action.
Eventually we come up with a further harmonic datum:
\[ p_{22} = 
  \langle x_0^4-6x_0^2x_1^2 +x_1^4, f_t \rangle ^2
+ 16 \langle x_0x_1(x_0^2-x_1^2), f_t \rangle ^2 \, . \]
Using this harmonic datum gives by Proposition \ref{data-theta}:

\begin{corollary}\label{theta2}
For a integral lattice $\Lambda \subset \edop^2$ the modular
form
\[\Theta_{22,\Lambda} =
\Theta_{x_0^4-6x_0^2x_1^2+x_1^4,\Lambda}^2
+16\Theta_{x_0x_1(x_0^2-x_1^2),\Lambda}^2
\]
of weight 10 is invariant under the $\Orth(2)$ action on the embedding $\Lambda \to
\edop^2$. 
The $q$-expansion of $\Theta_{22,\Lambda} = \sum_{n \geq 1} a_nq^n$ is
given by
\[a_n=
\sum_{\tiny \begin{array}{c}
(\gamma,\delta)  \in \Lambda \times \Lambda\\
||\gamma||^2+||\delta ||^2=n\\
\end{array}} \cos(4\measuredangle
(\gamma,\delta ))||\gamma||^4||\delta ||^4. \]
\end{corollary}

\begin{proof}
It is enough to derive from the
expression $\Theta_{22,\Lambda} =
\Theta_{x_0^4-6x_0^2x_1^2+x_1^4,\Lambda}^2
+16\Theta_{x_0x_1(x_0^2-x_1^2),\Lambda}^2$ the stated $q$-expansion
which is an elementary calculation.
However, since this is a special case of Theorem \ref{theta3}, we
omit the proof.
\end{proof}

\begin{theorem}\label{theta3}
For a positive integer $n\in \zdop$ we take the two harmonic polynomials
$h_1(x,y) = \Re((x+iy)^{2n})$, and $h_2(x,y)=\Im((x+iy)^{2n})$ of degree
$2n$.
For any integral lattice $\Lambda \subset \edop^2$ of level $N$
there is a $\Orth(2)$-invariant modular form 
\[ \Theta_{nn, \Lambda}(\tau) := \Theta_{h_1,\Lambda}^2(\tau) +
\Theta_{h_2,\Lambda}^2(\tau) \]
of level $N$, and weight $2+4n$. Its $q$-expansion
$\Theta_{nn, \Lambda}(\tau) = \sum_{n \geq 0} a_nq^n$ is given by
\[a_n=
\sum_{\tiny \begin{array}{c}
(\gamma,\delta)  \in \Lambda \times \Lambda\\
||\gamma||^2+||\delta ||^2=n\\
\end{array}} \cos(2n\measuredangle
(\gamma,\delta ))||\gamma||^{2n}||\delta ||^{2n}. \]
Furthermore, $\Theta_{nn, \Lambda} \in (4q^{2k}) \subset \zdop[[q]]$ where
$k$ is the first minimal length of the lattice.
\end{theorem}
\begin{proof}
We compute the $q$-expansion of $\Theta_{nn,\Lambda}$.
\[ \begin{array}{rcl}
\Theta_{nn,\Lambda}(\tau)
& =&
\sum\limits_{(\gamma,\delta ) \in
\Lambda \times \Lambda}\left(
h_1(\gamma)h_1(\delta)+h_2(\gamma)h_2(\delta)
\right)q^{||\gamma||^2+||\delta||^2}\\
& =&
\sum\limits_{(\gamma,\delta ) \in
\Lambda \times \Lambda}\left(
\Re(\gamma^{2n})\Re(\delta^{2n})+\Im(\gamma^{2n})\Im(\delta^{2n})
\right)q^{||\gamma||^2+||\delta||^2}\\
& =&
\sum\limits_{(\gamma,\delta ) \in
\Lambda \times \Lambda}\left(
\Re(\gamma^{2n} \bar \delta^{2n})
\right)q^{||\gamma||^2+||\delta||^2}\\
& =&
\sum\limits_{(\gamma,\delta ) \in
\Lambda \times \Lambda}\left(
\Re\left( \frac{\gamma^n \bar \delta^n}{\bar \gamma^n \delta^n}
\right) \gamma^n\bar\gamma^n\delta^n\bar\delta^n
\right)q^{||\gamma||^2+||\delta||^2}\\
& =&
\sum\limits_{(\gamma,\delta ) \in
\Lambda \times \Lambda}\left(
\cos(2n\measuredangle
(\gamma,\delta ))||\gamma||^{2n}||\delta ||^{2n}
\right)q^{||\gamma||^2+||\delta||^2}\\
\end{array}\]
We deduce from this form of the $q$-expansion that
$a_0=a_1=\dots=a_{2k-1}=0$, and that
$\Theta_{nn,\Lambda}$ does not depend from the embedding $\Lambda
\subset \edop^2$. It remains to demonstrate that $a_n \in 4\,\zdop$.

To show this, we use from the above deduction that $a_n =
\sum_{(\gamma,\delta) \in \Lambda \times \Lambda}
\Re(\gamma^{2n}\bar \delta^{2n})$ 
where we sum over all pairs $(\gamma,\delta)$ with
$||\gamma||^2+||\delta||^2=n$.
The symmetric bilinear form $\psi:\Lambda \times \Lambda \to \rdop$
with $\psi(\gamma,\delta):= \Re(\gamma \bar \delta)$ satisfies
$\psi(\gamma,\gamma) \in \zdop$ because $\Lambda $ is integral. Thus,
$\psi(\gamma,\delta) \in \frac{1}{2} \zdop$. We deduce that for
$\gamma,\delta \in \Lambda$ the element $x:=\gamma\bar\delta$ satisfies
a quadratic equation $x^2-b_1x+b_2=0$ with $b_1=2\Re(\gamma\bar\delta)
\in \zdop$, and $b_2= ||\gamma||^2||\delta||^2 \in \zdop$. Thus, $x$ is
an integer in a (at most) quadratic field extension $K/\qdop$. 
Consequently, $x^{2n}$ is an integer in this field $K$ which implies
$2\Re(\gamma^{2n}\bar\delta^{2n}) \in \zdop$.
If $\gamma \ne \pm \delta$, then the summand
$\Re(\gamma^{2n}\bar\delta^{2n})$ appears eight times in $a_n$, namely
from the pairs $(\pm \gamma, \pm \delta)$ and $(\pm \delta, \pm
\gamma)$. If $\gamma = \pm \delta$, then
$\Re(\gamma^{2n}\bar\delta^{2n})$ is an integer and appears four times.
\end{proof}

{\em Example 1. Some invariants for the lattice $\Lambda$
associated to $q(x,y)=x^2+y^2$}.
\[ \begin{array}{rcl}
\Theta_\Lambda(\tau) & = & 1 +4 q +4 q^{2} +4 q^{4} +8 q^{5} +4 q^{8}
+4 q^{9} +8 q^{10} +8 q^{13} +4 q^{16} +8 q^{17} + \dots\\
 \Theta_{22,\Lambda}(\tau)  &= &16(q^{2} -8 q^{3} +16 q^{4} +32 q^{5}
-156 q^{6} +112 q^{7} +256 q^{8} -576 q^{9} +\dots) \\
 \Theta_{44,\Lambda}(\tau)  &=& 16(q^{2} +32 q^{3} +256 q^{4} +512 q^{5}
+6084 q^{6} -33728 q^{7} +65536 q^{8} + \dots) \\
 \Theta_{11,\Lambda}(\tau)  &=&0 \,\, = \,\, \Theta_{33,\Lambda}(\tau)\\
\end{array}\]

{\em Example 2. Some invariants for the lattice $\Lambda$
associated to $q(x,y)=x^2+xy+y^2$}.
\[ \begin{array}{rcl}
\Theta_\Lambda(\tau)  &=& 1 +6 q +6 q^{3} +6 q^{4} +12 q^{7} +6 q^{9}
+6 q^{12} +12 q^{13} +6 q^{16} +12 q^{19} + \dots \\
\Theta_{33,\Lambda}(\tau) & =& 36(q^{2} -54 q^{4} +128 q^{5} +729 q^{6}
-3456 q^{7} +3524 q^{8} +16902 q^{10} + \dots )\\
\Theta_{11,\Lambda}(\tau) &= &0 \,\, = \,\, \Theta_{22,\Lambda}(\tau)
\,\, = \,\, \Theta_{44,\Lambda}(\tau)\\
\end{array}\]

\section{The harmonic datum $p_{11}$ for lattices in
$\edop^n$}\label{theta11-n}
\subsection{An explicit formula for $c_{11}$ and definition of $p_{11}$}
Let $\Lambda \subset \edop^n$ be a lattice. As explained in \ref{heat}
the heat flux of the distribution $\Lambda$ is given by
$\omega_{\Lambda,t} = f_t dx_1 \wedge \dots dx_n$ with
\[ f_t(x) = (4\pi t)^\frac{-n}{2}\sum_{\gamma \in \Lambda} \exp \left(
\frac{-||x-\gamma||^2}{4t} \right) .\]
The value of $f_t(0)$ is the first lattice invariant.
To derive the lattice invariant $c_{11}$ we need the second Taylor
coefficient of $f_t$ in zero.
We write it in the following manner:
\[ f_{t,1} = \sum_{i=1}^n \frac{a_i}{2}x_i^2 +
\sum_{1 \leq i < j \leq n} b_{ij}x_ix_j \]
where $a_i = \frac{\dee^2}{\dee x_i^2}f_t|_{x=0}$ and
$b_{ij}=\frac{\dee}{\dee x_i}\frac{\dee}{\dee x_j}f_t|_{x=0}$.
The invariant $c_{11}$ is defined to be the integral $c_{11}
:=\int_{S^{n-1}} f_{t,1}^2 d \bar \mu$.

As a preparation we compute $c_1:= \int_{S^{n-1}} f_{t,1} d \bar \mu$.
Using the integral formulas $\int_{S^{n-1}}x_i^2 d\bar \mu =\frac{1}{n}$
and $\int_{S^{n-1}}x_ix_j d\bar \mu =0$ (cf. Corollary \ref{int-pol}) we
see that
\[ c_1 = \frac{1}{2n} (\sum_{i=1}^n a_i) = \frac{-1}{2n} \Delta
f_t|_{x=0} = \frac{1}{2n} \frac{\dee f_t}{\dee t}|_{x=0} .\]

We need the following integrals from Corollary \ref{int-pol}:
\[ \int_{S^{n-1}}x_i^4 d\bar \mu = \frac{3}{n(n+2)} \, ,\qquad
\int_{S^{n-1}}x_i^2x_j^2 d\bar\mu =\frac{1}{n(n+2)} \, . \]
Having in mind (Corollary \ref{int-pol}) that the integrals of
homogeneous monomials over $S^{n-1}$ vanish if an odd exponent occurs we
find the following:
\[\begin{array}{rcl}
4n(n+2) c_{11} & =& \displaystyle
4\sum_{1 \leq i < j \leq n} b_{ij}^2 + 3 \sum_{i=1}^n a_i^2
+ 2\sum_{1 \leq i < j \leq n} a_i a_j\\
\\
&=& \displaystyle
4\sum_{1 \leq i < j \leq n} b_{ij}^2 + 2 \sum_{i=1}^n a_i^2
+ \sum_{i=1}^n\sum_{j=1}^n a_i a_j .
\end{array}\]
We want to express $c_{11}$ in terms of differential operators which
correspond either to harmonic polynomials or to
$r:= a_1+ a_2 + \dots +a_n $. This $r$ corresponds to $-\Delta$, as
$r= -\Delta(f)|_{x=0}$.
So we introduce the harmonic polynomials $h_i= na_i -r$.
Using the obvious equality $\sum_{i=1}^n h_i =0$ we obtain:
\[\begin{array}{rcl}
4n^3(n+2) c_{11} & =& \displaystyle
4n^2\sum_{1 \leq i < j \leq n} b_{ij}^2 + 2 \sum_{i=1}^n (na_i)^2
+ \sum_{i=1}^n\sum_{j=1}^n (n a_i)(n a_j) \\
& =& \displaystyle
4n^2\sum_{1 \leq i < j \leq n} b_{ij}^2 +2 \sum_{i=1}^n (h_i+r)^2
+ \sum_{i=1}^n\sum_{j=1}^n (h_i+r)(h_j+r) \\
& =& \displaystyle
\left (4n^2\sum_{1 \leq i < j \leq n} b_{ij}^2 +2 \sum_{i=1}^n h_i^2
\right) 
+\left( (n^2+2n)r^2 \right) \,.\\
\end{array}\]

So the invariant $c_{11}$ decomposes as $c_{11}=
\frac{1}{2n^3(n+2)}p_{11} + \frac{r^2}{4n^2}$.
The second summand corresponds to $c_1^2 =
\left(\frac{1}{2n}\frac{\dee}{\dee r}f_t|_{x=0}\right)^2$.
Hence, the first summand is also invariant.
It yields a harmonic datum:
\[ p_{11} = 
2n^2\sum_{1 \leq i < j \leq n} b_{ij}^2 +
\sum_{i=1}^n h_i^2 = 2n^2\sum_{1 \leq i < j \leq n} \langle x_ix_j,f_t
\rangle^2 + \sum_{i=1}^n \langle h_i, f_t \rangle^2
.\]
{\em Remark.}
Note that for $n=2$ we have $h_1=-h_2$, so $h_1^2=h_2^2$. Thus, in this
case we have $p_{11} = 2\left( 4 \left(\frac{\dee}{\dee
x_1}\frac{\dee}{\dee x_2}f_t|_{x=0}\right)^2 + \left((\frac{\dee^2}{\dee x_1^2} - \frac{\dee}{\dee
x_2^2} ) f_t|_{x=0}\right)^2 \right) $.  Up to the factor 2 this is
our harmonic system from \ref{p112}.\\

Applying Proposition \ref{data-theta} we obtain the next result:

\begin{theorem}\label{theta11}
Let $\Lambda \subset \edop^n$ be an integral lattice,
of level $N$.
The modular form
\[ \Theta_{11,\Lambda}(\tau) = 2n^2\left( \sum_{1 \leq i < j \leq n}
\Theta_{x_ix_j,\Lambda}^2(\tau) \right) + \sum_{i=1}^n
\Theta_{nx_i^2-\sum_{j=1}^nx_j^2,\Lambda}^2(\tau) \]
is independent from the embedding $\Lambda \to \edop^n$.
$\Theta_{11,\Lambda}$ is a cusp form of weight $n+4$ of level $N$.
Its $q$-expansion is given by 
\[\Theta_{11,\Lambda}(\tau) = \sum_{m \geq 0}a_m q^m \quad with \quad
a_m=
n^2  \! \! \! \! \! \! \! \! \! \! \! \!
\sum_{\tiny \begin{array}{c}
(\gamma,\delta) \in \Lambda \times \Lambda\\
||\gamma||^2+||\delta ||^2=m\\
\end{array}} \! \! \! \! \! \! \! \! \! \! \left(\cos^2(\measuredangle
(\gamma,\delta ))-\frac{1}{n}\right)||\gamma||^2||\delta ||^2. \]
We have $\Theta_{11,\Lambda} \in (2nq^{2k}) \subset \zdop[[q]]$
with $k$ the first minimal length of the lattice $\Lambda$.
If $n$ is even, then $\Theta_{11,\Lambda} \in (4nq^{2k}) \subset
\zdop[[q]]$.
\end{theorem}

\begin{proof} On the one hand, when considering
the $q$-expansion it is obvious that
$\Theta_{11,\Lambda}$ is invariant under the $\Orth(n)$ action.
On the other hand, being a sum of squares of modular forms of weight
$\frac{n}{2}+2$ forces $\Theta_{11,\Lambda}$ to be a modular form of
weight $n+4$. Thus, even though $\Theta_{11,\Lambda}$ comes from a
harmonic datum, all we have to do is to show that the modular form has
the stated $q$-expansion.
\[ \begin{array}{rcl}
\Theta_{11,\Lambda}(\tau)
& =& \! \! \! \! \! \! \! \!
\sum\limits_{(\gamma,\delta) \in
\Lambda \times \Lambda}\left( 2n^2 \sum\limits_{1 \leq i < j \leq n}
\gamma_i\gamma_j\delta_i\delta_j  + \sum\limits_{i=1}^n
(n\gamma_i^2-||\gamma||^2)(n\delta_i^{2}-||\delta ||^2)
\right)q^{||\gamma||^2+||\delta||^2}\\
& =& \! \! \! \! \! \! \! \!
\sum\limits_{(\gamma,\delta) \in
\Lambda \times \Lambda}\left( n^2 \sum\limits_{i=1}^n
\sum\limits_{j=1}^n
\gamma_i\gamma_j\delta_i\delta_j   -n ||\gamma||^2||
\delta||^2 \right)q^{||\gamma||^2+||\delta||^2}\\
& =& \! \! \! \! \! \! \! \!
\sum\limits_{(\gamma,\delta) \in
\Lambda \times \Lambda}\left( n^2 
\langle \gamma, \delta \rangle^2  -n ||\gamma||^2||
\delta||^2 \right)q^{||\gamma||^2+||\delta|^2}\\
\end{array}\]
Now the definition of the cosine gives the formula for the
$q$-expansion. From this formula we conclude that $a_0=a_1=\dots
=a_{2k-1}=0$.
In order to prove that $a_m \in 2n \zdop$ we consider the sum
\[\frac{1}{n}a_m=
\sum_{\tiny \begin{array}{c}
(\gamma,\delta) \in \Lambda \times \Lambda\\
||\gamma||^2+||\delta ||^2=m\\
\end{array}}
\left( n\langle \gamma ,\delta
\rangle^2-||\gamma||^2||\delta||^2 \right)
\]
Since $\Lambda$ is integral we have $n\langle \gamma ,\delta
\rangle^2 \in \frac{1}{4} \zdop$ (respectively in $\frac{1}{2}\zdop$
when $n$ is even). If $\gamma \ne \pm \delta$, then the eight pairs
$(\pm \gamma, \pm \delta)$ and $(\pm \delta, \pm \gamma)$ in $\Lambda
\times \Lambda$ give the same contribution
$\left( n\langle \gamma ,\delta \rangle^2-||\gamma||^2||\delta||^2
\right)$
 to $a_m$.
If $\gamma \ne \pm \delta$, then $\langle \gamma ,\delta
\rangle$ is an integer and the integer summand
$\left( n\langle \gamma ,\delta
\rangle^2-||\gamma||^2||\delta||^2 \right)$ appears four times.
\end{proof}

\subsection{Example: Computing $\Theta_{11,\Lambda}$ for two isospectral
lattices in dimension four}

We consider the two integral lattices $\Lambda_1$ and $\Lambda_2$ in
$\edop^4$ which were investigated by A.~Schiemann in \cite{Sch}.
The Gram matrices of these lattices are
given by
\[ A_1 = \frac{1}{2} \left( \begin {array}{cccc}
4&2&0&1\\\noalign{\medskip}2&8&3&1\\\noalign{\medskip}0&3&10&5\\\noalign{\medskip}1&1&5&10\end
{array} \right) \qquad
A_2 = \frac{1}{2}\left( \begin {array}{cccc}
4&0&1&1\\\noalign{\medskip}0&8&1&-4\\\noalign{\medskip}1&1&8&2\\\noalign{\medskip}1&-4&2&10\end
{array} \right) \]
Schiemann showed in his article that these two lattices are not
isometric even though they give the same theta function $\Theta(\tau)
=\Theta_{\Lambda_1}(\tau) = \Theta_{\Lambda_2}(\tau)$ which he
determined to be
\[ \Theta(\tau)
=  1 +2 q^2 +4 q^4 +6 q^5 +10 q^6 +6 q^7 +12 q^8 +6 q^9 +6 q^{10}
+8 q^{11} +10 q^{12} +8 q^{13} +10 q^{14} +22 q^{15} + \dots\]

\begin{proposition}\label{exam-schie}
We have an inequality
$\Theta_{11,\Lambda_1} \ne \Theta_{11,\Lambda_2}$.
Thus, the modular forms $\Theta_{11,\Lambda_1}$ and
$\Theta_{11,\Lambda_2}$ of level 1729 and weight eight distinguish the two
isospectral lattices $\Lambda_1$ and $\Lambda_2$.
\end{proposition}
\begin{proof}
We need embeddings of $\Lambda_i \to \edop^4$. We choose a
decomposition $A_i = S_i^t \cdot S_i$ with $S_i$ upper triangular.
We obtain
\[ S_1 = \frac{1}{\sqrt{2}}
\left( \begin {array}{cccc} 2&1&0&\frac{1}{2}\\
\noalign{\medskip}0&\sqrt {7}&\frac{3\sqrt{7}}{7}&\frac{\sqrt {7}}{14}\\
\noalign{\medskip}0&0&\frac{\sqrt {427}}{7}&{
\frac {67 \sqrt {427}}{854}}\\\noalign{\medskip}0&0&0&{\frac
{\sqrt {105469}}{122}}\end {array} \right) \]
The column vectors $\{ \gamma_i \}_{i=1,2,3,4}$ of the matrix $S_1$
generate the lattice $\Lambda_1$.
From the theta function we see that there are two vectors of norm
$\sqrt{2}$ ($\pm \gamma_1$), and four vectors of norm two ($\pm
\gamma_2$, and $\pm(\gamma_1-\gamma_2)$). Fortunately, these six lattice
vectors are enough. We compute now in the ring $A = \rdop[[q]]$.
The six spherical theta functions $\Theta_{\Lambda_1,x_ix_j}$ for $i<j$
have zero coefficient at $q^2$. Since there are no lattice vectors of
norm $\sqrt{3}$ we find $\Theta_{\Lambda_1,x_ix_j} \in (q^4)$. Thus,
we deduce $\Theta_{\Lambda_1,x_ix_j}^2 \in (q^8)$.

Considering the lattice vectors of norm at most 2, we see
\[ \begin{array}{rcl}
\Theta_{\Lambda_1,4x_1^2-(x_1^2+x_2^2+x_3^2+x_4^2)} & = & 12q^2 -8q^4 +
\dots \\
\Theta_{\Lambda_1,4x_2^2-(x_1^2+x_2^2+x_3^2+x_4^2)} & = & -4q^2 +40q^4 +
\dots \\
\Theta_{\Lambda_1,4x_3^2-(x_1^2+x_2^2+x_3^2+x_4^2)} & = & -4q^2 -16q^4 +
\dots \\
\Theta_{\Lambda_1,4x_4^2-(x_1^2+x_2^2+x_3^2+x_4^2)} & = & -4q^2 -16q^4 +
\dots \\
\end{array}\]
Taking the squares of these equations we obtain
\[ \begin{array}{rcl}
\Theta^2_{\Lambda_1,4x_1^2-(x_1^2+x_2^2+x_3^2+x_4^2)} & = & 144q^4 -192q^6 +
\dots \\
\Theta^2_{\Lambda_1,4x_2^2-(x_1^2+x_2^2+x_3^2+x_4^2)} & = & 16q^4 -320q^6 +
\dots \\
\Theta^2_{\Lambda_1,4x_3^2-(x_1^2+x_2^2+x_3^2+x_4^2)} & = & 16q^4 +128q^6 +
\dots \\
\Theta^2_{\Lambda_1,4x_4^2-(x_1^2+x_2^2+x_3^2+x_4^2)} & = & 16q^4 +128q^6 +
\dots \, .\\
\end{array}\]
Summing up these four squares we obtain that
\[ \Theta_{11,\Lambda_1}(\tau) = 192q^4-256q^6 + \dots \,.\]
We repeat this construction now with the lattice $\Lambda_2$.
In this case we find
\[ S_2 = \frac{1}{\sqrt{2}} 
 \left( \begin {array}{cccc} 2&0&\frac{1}{2}&\frac{1}{2}\\
\noalign{\medskip}0&2\,\sqrt
{2}&\frac{\sqrt {2}}{4}&-\sqrt {2}\\
\noalign{\medskip}0&0&\frac{\sqrt {122}}{4}&{
\frac {9 \sqrt {122}}{122}}\\
\noalign{\medskip}0&0&0&{\frac
{\sqrt {105469}}{122}}\end {array} \right) \]
The vectors of norm $\sqrt{2}$ are $\pm \gamma_1$, whereas the vectors
of norm two are the four vectors $\pm \gamma_2$, and $\pm \gamma_3$.
Again, we have $\Theta_{\Lambda_2,x_ix_j}^2 \in (q^8)$ for $i<j$.
As before we compute
\[ \begin{array}{rcl}
\Theta_{\Lambda_2,4x_1^2-(x_1^2+x_2^2+x_3^2+x_4^2)} & = & 12q^2 -15q^4 +
\dots \\
\Theta_{\Lambda_2,4x_2^2-(x_1^2+x_2^2+x_3^2+x_4^2)} & = & -4q^2
+\frac{33}{2}q^4 +
\dots \\
\Theta_{\Lambda_2,4x_3^2-(x_1^2+x_2^2+x_3^2+x_4^2)} & = & -4q^2
+\frac{29}{2}q^4 +
\dots \\
\Theta_{\Lambda_2,4x_4^2-(x_1^2+x_2^2+x_3^2+x_4^2)} & = & -4q^2 -16q^4 +
\dots  \,. \\
\end{array}\]
Eventually, we obtain $\Theta_{11,\Lambda_2}(\tau)= 192q^4-480q^6+ \dots
\ne \Theta_{11,\Lambda_1}(\tau)$.
\end{proof}

{\em Remark.}
To see that the two modular forms are different it is enough to compute
up to the first Fourier coefficient which is different.
Anyway, using a computer we can give more coefficients:
$$\begin{array}{rcl}
\Theta_{11,\Lambda_1}(\tau) &=&
192 q^4 -256 q^6 -896 q^7 +1120 q^8 -2848 q^9 +3024 q^{10} -2112 q^{11}
+\\
&& \qquad \qquad \qquad \qquad +13536 q^{12} -4064 q^{13} -16272 q^{14} -4544 q^{15} + \dots\\
\Theta_{11,\Lambda_2}(\tau) &=&
192 q^4 -480 q^6 -608 q^7 +736 q^8 -1312 q^9 +3216 q^{10} +1056
q^{11}-\\
&&\qquad \qquad \qquad \qquad -2048 q^{12} -2624 q^{13} +2896 q^{14} -12288 q^{15} +\dots\\
\end{array}$$

\appendix
\section{Integrating polynomials on spheres}\label{int-pol-sn}
Let $f \in \rdop[x_0,x_1,\dots,x_n]$ be a polynomial. We need the
integral $\int_{S^n} f d\mu$ for further computations. Here $d\mu$
denotes the standard $\Orth(n+1)$ invariant measure on $S^n$.
Since we easily can decompose
$f$ into its homogeneous components, it is enough to consider
homogeneous polynomials $f$. Here we have the following result:
\begin{proposition}\label{integral}
Let $f \in \rdop[x_0,x_1,\dots,x_n]$ be an homogeneous polynomial of
degree $d$. If $d$ is odd, we have that $\int_{S^n} f d \mu =0$.
If $d =2k$ is even, $\Delta^k f$ is a real number, and we
have
\[ \int_{S^n} f d \mu = c_d \Delta^k f \quad \mbox{with} \quad 
c_d:=\alpha_d\cdot 
\frac{(n+1) \pi^{\frac{n+1}{2}}}{\Gamma(\frac{n+3}{2})}
\mbox{ , and }  \alpha_{2k}:= \frac{1}{(-2)^k k! \prod_{m=1}^k
(n+2m-1)} \,.\]
\end{proposition}
\begin{proof}
We start with the following observation (see Exercise 33 on page 550 in
\cite{Lang}): $f$ can be uniquely decomposed 
\[ f = \sum_{l=0}^{\lfloor \frac{d}{2} \rfloor} r^l h_{d-2l} \]
where $r=\sum_{i=0}^n x_i^2$ and the $h_i$ are homogeneous
harmonic polynomials of degree $i$, that is $\Delta h_i =0$.
The mean value principle and the fact that $r$ is constantly one on $S^n$
together with the fact that the volume of $S^n$ is 
$\frac{(n+1) \pi^{\frac{n+1}{2}}}{\Gamma(\frac{n+3}{2})}$ yield
\[ \int _{S^n} f = \int _{S^n} \sum_{l=0}^{d/2} h_{d-2l}=
\frac{(n+1) \pi^{\frac{n+1}{2}}}{\Gamma(\frac{n+3}{2})}
\sum_{l=0}^{d/2} h_{d-2l}(0)  \]
Now the value of an homogeneous polynomial at zero is zero, unless the
degree of the polynomial is zero. Thus, we deduce that 
$\int_{S^n} f d \mu = \frac{(n+1)
\pi^{\frac{n+1}{2}}}{\Gamma(\frac{n+3}{2})}h_0$. This, shows that
integration of homogeneous polynomials of odd degree gives zero. Next,
we show (by simple and straightforward differentiation) the following equality 
for any homogeneous harmonic polynomial $h$:
\[ \Delta( r^m h) = (-2m) (n+2m+2 \deg(h) -1)r^{m-1}h \,.\]

In consequence we have $\Delta^k r^k = (-2)^k k! \prod_{m=1}^k (n+2m-1)
= \frac{1}{\alpha_{2k}}$, and $\Delta^kr^m h=0$ for $m<k$ and $h$ harmonic.
We deduce that $\Delta^k f$ yields $\frac{h_0}{\alpha_{2k}}$, which
finishes the proof.
\end{proof}

We need concrete formulas for the normalized measure $d \bar \mu =
\frac{\Gamma(\frac{n+3}{2})}{(n+1)\pi^\frac{n+1}{2}} d \mu$. This measure
is $\Orth(n)$-invariant and has the property that $\int_{S^n} d \bar \mu
=1$. As an application of the above Proposition \ref{integral} we
obtain:

\begin{corollary}\label{int-pol}
For the monomial $f = x_0^{i_0}x_1^{i_1} \cdots x_n^{i_n}$ of degree
$d=i_0+i_1+\dots  + i_n$ we have the equality
\[ \int _{S^n} f d \bar \mu = \left\{ \begin{array}{ll}
0 & \mbox{at least one of the }i_k\mbox{ is odd.}\\
\\
\displaystyle
\left( \prod_{k=0}^n\frac{i_k!}{\frac{i_k}{2}!} \right)
\left( \prod_{m=1}^{d/2} \frac{1}{2(n+2m-1)} \right) \quad &
\mbox{all the }i_k\mbox{ are even.}\\
\end{array} \right.\]

\end{corollary}


\begin{thebibliography}{1}
\bibitem{BGV} N.~Berline, E.~Getzler, and M.~Vergne,
{\em Heat Kernel and Dirac Operators},
Grundlehren der math.~Wiss.~{\bf 298}, Springer, Berlin, 1996.
\bibitem{Hei} G. Hein,
{\em Computing Green Currents via the Heat Kernel},
J. reine angew. Math. {\bf 540} (2001) 87-104.
\bibitem{Lang} S.~Lang,
Algebra, revised third ed., Springer, New York, 2002.
\bibitem{Sch} A.~Schiemann, {\em Ein Beispiel positiv definiter
quadratischer Formen der Dimension 4 mit gleichen Darstellungszahlen},
Arch.~Math.~{\bf 54} (1990) 372--375.
\bibitem{Zag} D.~Zagier,
{\em Introduction to modular forms} in:
M.~Waldschmidt, P.~Moussa, J.~M.~Luck and C.~Itzykson (eds.)
 From number theory to physics (Les
Houches, 1989), 238--291, Springer, Berlin, 1992. 
\end{thebibliography}
\end{document}